\documentclass[12pt, a4paper]{article}

\usepackage[utf8]{inputenc}
\usepackage[T1]{fontenc}
\usepackage{geometry}
\usepackage{amsmath, amssymb, amsfonts, amsthm}
\usepackage{mathtools}

\usepackage{graphicx}
\usepackage{float}
\usepackage{booktabs}
\usepackage{caption}
\usepackage{subcaption}
\usepackage{tikz}
\usepackage{pgfplots}
\pgfplotsset{compat=1.18}
\usetikzlibrary{positioning, arrows.meta, calc, shapes.geometric}

\usepackage[ruled,vlined]{algorithm2e}

\usepackage[authoryear]{natbib}

\usepackage{authblk}

\usepackage{xcolor}
\usepackage[colorlinks=true, linkcolor=blue, citecolor=blue, urlcolor=blue]{hyperref}

\theoremstyle{plain}
\newtheorem{theorem}{Theorem}[section]
\newtheorem{lemma}[theorem]{Lemma}
\newtheorem{corollary}[theorem]{Corollary}
\newtheorem{proposition}[theorem]{Proposition}

\newtheorem{Assumption}[theorem]{Assumption}

\theoremstyle{definition}

\theoremstyle{remark}
\newtheorem{remark}[theorem]{Remark}

\title{\textbf{Optimal Control of an Epidemic with Intervention Design}}
\author[1]{Behrooz Moosavi Ramezanzadeh\thanks{Email: \texttt{behroozmoosavi@pitt.edu}}
\thanks{The author is grateful to two anonymous referees for their careful reading and constructive comments, which significantly improved the manuscript.}}
\affil[1]{Ph.D. Student in Economics \\ University of Pittsburgh, Pittsburgh, PA, USA}
\date{\today}

\begin{document}

\maketitle

\begin{abstract}
This paper investigates the optimal control of an epidemic governed by a SEIR model with an operational delay in vaccination. We address the mathematical challenge of imposing hard healthcare capacity constraints (e.g., ICU limits) over an infinite time horizon. To rigorously bridge the gap between theoretical constraints and numerical tractability, we employ a variational framework based on Moreau--Yosida regularization and establish the connection between finite- and infinite-horizon solutions via $\Gamma$-convergence. The necessary conditions for optimality are derived using the Pontryagin Maximum Principle, allowing for the characterization of boundary-maintenance arcs where the optimal strategy maintains the infection level precisely at the capacity boundary. Numerical simulations illustrate these theoretical findings, quantifying the shadow prices of infection and costs associated with intervention delays.
\end{abstract}


\section{Introduction}

The management of infectious disease outbreaks—from the historical 1918 Spanish Flu to the contemporary global crisis of COVID-19—represents a fundamental challenge for public health stability and economic prosperity. These catastrophic events have demonstrated that epidemic control is not merely a biological problem but a complex intertemporal optimization task. Public health authorities must design intervention strategies that balance the immediate goal of minimizing infection peaks against the long-term socio-economic costs of control, all while navigating rigid operational constraints such as finite vaccine supplies, implementation delays, and the inelastic capacity of healthcare infrastructure \cite{barry-2004, ferguson-2020}.

To address these trade-offs, mathematical modeling has historically relied on compartmental frameworks, most notably the Susceptible-Exposed-Infected-Recovered (SEIR) model \cite{anderson-1991, kermack-1927}. While descriptive forecasting has provided essential insights, the research focus has increasingly shifted toward prescriptive optimal control theory. This approach seeks to determine the optimal timing and intensity of interventions to prevent healthcare system collapse under realistic resource limitations \cite{bussell-2019, behncke-2000}. However, a significant disconnect persists in the literature between theoretical optimality and operational reality. Many existing models treat healthcare capacity constraints either as "hard" boundaries that complicate numerical solutions or as "soft" penalties that lack rigorous mathematical justification in an infinite-horizon setting. Furthermore, the economic interpretation of control terms often conflates distinct mechanisms, such as distinguishing between therapeutic hospitalization and transmission-reducing suppression.

This paper addresses these limitations by establishing a rigorous variational framework that explicitly bridges the gap between theoretical state constraints and numerical policy design. Building on recent foundational work in infinite-horizon epidemic control \cite{freddi-2022, freddi-goreac-2024, della-rossa-2025}, we contribute to the literature in three primary ways.

First, we provide a mathematically well-posed SEIR model that incorporates realistic operational delays and proves essential stability properties. Unlike standard formulations that may violate biological invariance under delayed controls, we utilize Barbalat's Lemma to establish the asymptotic convergence of the infected population to zero, ensuring the long-term validity of the control problem.

Second, we resolve the ambiguity surrounding healthcare capacity constraints. We formulate the problem with a "hard" state constraint $i(t) \le I_{\max}$ to represent absolute ICU capacity limits. To facilitate numerical analysis without sacrificing theoretical rigor, we utilize Moreau--Yosida regularization to establish a formal link between this hard constraint and a penalized objective function \cite{DalMaso1993, RockafellarWets1998}. We prove that the solution to the penalized problem converges to the hard-constrained optimum as the penalty parameter increases, ensuring that our numerical results are theoretically consistent with strict infrastructure limits.

Third, we clarify the economic and biological interpretation of the control mechanisms. We distinguish between vaccination $u(t)$ and transmission suppression $h(t)$. Unlike interpretations that conflate suppression with hospitalization, we define $h(t)$ as \textbf{Non-Pharmaceutical Interventions (NPIs)} (e.g., social distancing, testing, isolation, lockdowns) that reduce the effective transmission rate $\beta$. Accordingly, we justify the structure of our cost functional based on economic principles: vaccination costs are modeled linearly proportional to administration intensity (through the term $c_V u(t)s(t)$), and suppression costs are modeled as prevalence-scaled resource burdens (through the term $c_H i(t)h(t)$), capturing the fact that enforcement and tracing resources are most intensively deployed when infections are widespread.

Finally, we employ $\Gamma$-convergence to rigorously justify the transition from finite-horizon numerical simulations to infinite-horizon theoretical stability. Extending the recent results of Freddi and Goreac \cite{freddi-goreac-2024} and Della Rossa and Freddi \cite{della-rossa-2025} to the SEIR context with delays, we show that the optimal strategies obtained over large finite horizons approximate the true infinite-horizon optimum. This analysis allows us to characterize boundary-maintenance arcs, where the system is optimally held at capacity—providing a meaningful characterization of the optimal strategy that descriptive models often miss.

\textbf{Plan of Paper:} The manuscript is organized as follows. Section 2 establishes the mathematical foundations of the controlled SEIR system, proving the existence and non-negativity of solutions. Section 3 formulates the dynamic optimization problem, introduces the variational bridge for state constraints via regularization, and derives optimality conditions using the Pontryagin Maximum Principle. Section 4 validates the framework through numerical simulations and sensitivity analyses. Finally, Section 5 offers concluding remarks.

\textbf{Notation:} Throughout this paper, $I=[0, \infty)$ denotes the infinite time horizon. We utilize the space $L^{\infty}(I)$ for essentially bounded measurable functions and $W^{1, \infty}(I; \mathbb{R}^k)$ for Lipschitz continuous state trajectories. The subscript ``$\mathrm{loc}$'' denotes the local variant: $f \in W^{1,\infty}_{\mathrm{loc}}(I; \mathbb{R}^k)$ if $f \in W^{1,\infty}([0,T]; \mathbb{R}^k)$ for every $T > 0$. The notation $u_n \rightharpoonup^\ast u$ denotes weak-* convergence in the control space. For a thorough mathematical treatment of these functional spaces, we refer the reader to \cite{fonseca-2006} and \cite{leoni-2017}.
\section{Analysis of Controlled System}
\label{sec:basic_results}

This section establishes the mathematical foundations of the model we intend to investigate. To achieve this, we focus on the SEIR model framework (see \cite{anderson-1991} for details). The controlled dynamics are defined by the following system of differential equations:

\begin{equation}\label{sdeq:1}
\begin{aligned}
& \frac{d s}{d t}=-(\beta-h(t)) s(t) i(t)-u(t) s(t), \\
& \frac{d e}{d t}=(\beta-h(t)) s(t) i(t)-\sigma e(t), \\
& \frac{d i}{d t}=\sigma e(t)-\gamma i(t), \\
& \frac{d r}{d t}=\gamma i(t)+u(t) s(t),
\end{aligned}
\end{equation}
subject to the population conservation identity:
\[
s(t)+e(t)+i(t)+r(t)=1.
\]

The control variables are constrained as follows:
\begin{itemize}
    \item The vaccination rate $u(t)$ satisfies $0 \leq u(t) \leq u_{\max }$ and is subject to a delay such that $u(t)=0$ for $t<t_{\text {delay }}$.
    \item The suppression effort $h(t)$ satisfies $0 \leq h(t) \leq h_{\max }$. We interpret $h(t)$ as a rate-reducing suppression instrument (representing the intensity of measures such as distancing, testing, or isolation) that subtracts from the baseline transmission rate $\beta$; it is not defined as a dimensionless share or fraction.
\end{itemize}
We will use the shorthand admissible control set
\[
\mathcal{U}_{\mathrm{ad}}
:=\left\{(u,h)\in L^\infty([0,\infty);[0,u_{\max}])\times L^\infty([0,\infty);[0,h_{\max}]) \;:\;
u(t)=0 \ \text{for a.e. } t<t_{\mathrm{delay}}\right\}.
\]
Furthermore, we define a treatment capacity threshold $I_{\max} \in (0, 1)$. In the subsequent optimization analysis, this threshold will serve either as a boundary for the set of feasible states or as a reference parameter for a penalty term in the cost functional.

Here, $\beta$ is the baseline transmission rate, and the effective transmission rate $(\beta-h(t))$ decreases as the intervention $h(t)$ increases.

\begin{figure}[ht]
\centering
\begin{tikzpicture}[scale=0.9, every node/.style={scale=0.9}, >=stealth, thick]
    \node[draw, circle, minimum size=1cm] (s) at (90:3cm) {s};
    \node[draw, circle, minimum size=1cm] (e) at (0:3cm) {e};
    \node[draw, circle, minimum size=1cm] (i) at (270:3cm) {i};
    \node[draw, circle, minimum size=1cm] (r) at (180:3cm) {r};
    
    \draw[->, bend left=20] (s) to node[above right, xshift=-5pt] {Infection: \((\beta-h(t))\,si\)} (e);
    \draw[->, bend left=20] (e) to node[right] {Incubation: \(\sigma\,e\)} (i);
    \draw[->, bend left=20] (i) to node[below] {Recovery: \(\gamma\,i\)} (r);
    \draw[->, bend left=20] (s) to node[left] {Vaccination: \(u(t)s\)} (r);
    
    \node[draw, rectangle, fill=yellow!30, rounded corners, inner sep=3pt] (h) at ($(s)!0.4!(e)+(-0.5,-1.5cm)$) {Suppression \(h(t)\)};
    \draw[->, dashed] (h) -- ($(s)!0.5!(e)$);
    
    \node[draw, rectangle, fill=green!30, rounded corners, inner sep=3pt] (u) at ($(s)!0.5!(r)+(-1,1cm)$) {Vaccination \(u(t)\)};
    \draw[->, dashed] (u) -- ($(s)!0.5!(r)$);
    
    \node[align=center, below=1cm of i] {Threshold: \\ \(i(t)\le I_{\max}\)};
    \node[align=center, above=1cm of s] {\footnotesize Vaccination starts at \\ \(t\ge t_{\text{delay}}\)};
\end{tikzpicture}
\caption{Circular diagram of the SEIR epidemic model. Nodes represent compartments, arrows represent transitions, and colored boxes indicate control interventions.}\label{figure1}
\end{figure}

To ensure the mathematical well-posedness of the system, we make the following formal assumptions (see Hale \cite{hale1980ordinary} or Walter \cite{walter-1998}).

\begin{Assumption}[Model Parameters and Admissibility] \label{assump:bounds}
We assume the following conditions hold throughout the analysis:
\begin{enumerate}
    \item \textbf{Transmission and Suppression:} The baseline transmission rate $\beta$ and the maximum suppression effectiveness $h_{\max}$ satisfy $0 \le h_{\max} < \beta$. This ensures the effective transmission rate is strictly positive.
    \item \textbf{Control Regularity:} The controls satisfy $(u, h) \in L^{\infty}([0, \infty); [0, u_{\max}]) \times L^{\infty}([0, \infty); [0, h_{\max}])$.
    \item \textbf{Initial Conditions:} The initial populations satisfy $s_0 > 0$, $i_0 > 0$, $e_0 \ge 0$, $r_0 \ge 0$, and sum to 1.
\end{enumerate}
\end{Assumption}
\begin{theorem}[Existence, Uniqueness, and Non-negativity]\label{diffexist}
Let $I=[0,+\infty)$. Under Assumption \ref{assump:bounds}, there exists a unique global solution
$(s,e,i,r)\in W^{1,\infty}_{\mathrm{loc}}(I;\mathbb{R}^4)$ to \eqref{sdeq:1}. Moreover, the trajectory
remains in the compact simplex $\{(s,e,i,r)\in[0,1]^4:\ s+e+i+r=1\}$, so the vector field is uniformly
bounded; in particular there exists $K>0$ such that $\|(s,e,i,r)'\|_{L^\infty(I)} = \|(\dot{s},\dot{e},\dot{i},\dot{r})\|_{L^\infty(I)} \le K$.
For all $t\ge 0$,
\[
s(t)>0,\qquad e(t)\ge 0,\qquad i(t)\ge 0,\qquad r(t)\ge 0,
\]
and $i(t)>0$ for every finite $t\ge 0$ whenever $i_0>0$.
\end{theorem}

\begin{proof}
Consider the reduced system for $(s,e,i)$ on the set
$\mathbb{T}=\{(s,e,i)\in[0,1]^3:\ s+e+i\le 1\}$.
For each fixed control pair $(u,h)\in \mathcal{U}_{\mathrm{ad}}$, the right-hand side is measurable in $t$
and locally Lipschitz in $(s,e,i)$, hence the Carath\'eodory theorem yields a unique local absolutely
continuous solution.

Summing the four equations in \eqref{sdeq:1} gives $(s+e+i+r)'=0$, so $s(t)+e(t)+i(t)+r(t)\equiv 1$.
Consequently, as long as the components are nonnegative the trajectory stays in a compact subset of
$[0,1]^4$, where the vector field is bounded; thus the local solution extends globally to $[0,\infty)$.

\textbf{Step 1: strict positivity of $s$.}
Integrating the first equation gives
\[
s(t)=s_0 \exp\!\left\{-\int_0^t\big[(\beta-h(\xi))\,i(\xi)+u(\xi)\big]\,d\xi\right\}.
\]
Since $s_0>0$ and the integrand is bounded, $s(t)>0$ for all $t\ge 0$.

\textbf{Step 2: forward invariance of $e\ge 0$ and $i\ge 0$.}
Define the first-exit time
\[
\tau:=\inf\{t\ge 0:\ e(t)<0\ \text{or}\ i(t)<0\},
\qquad \inf\emptyset:=+\infty\quad 
\]

where the last convention means that if the set is empty—i.e., neither $e(t)$ nor $i(t)$ ever becomes negative—then $\tau=+\infty$ and the trajectory remains non-negative for all time.

By continuity and the initial conditions, $\tau>0$ and $e(t)\ge 0$, $i(t)\ge 0$ on $[0,\tau)$.
Assume for contradiction that $\tau<\infty$. Then at $t=\tau$ at least one component hits the boundary.
If $e(\tau)=0$ and $i(\tau)\ge 0$, then
\[
e'(\tau)=(\beta-h(\tau))\,s(\tau)\,i(\tau)\ge 0,
\]
so $e$ cannot cross from nonnegative to negative at $\tau$.
If $i(\tau)=0$ and $e(\tau)\ge 0$, then
\[
i'(\tau)=\sigma e(\tau)\ge 0,
\]
so $i$ cannot cross from nonnegative to negative at $\tau$.
Both cases contradict the definition of $\tau$, hence $\tau=+\infty$ and $e(t)\ge 0$, $i(t)\ge 0$
for all $t\ge 0$.

\textbf{Step 3: strict positivity of $i$ for finite times.}
Using variation of constants,
\[
i(t)=i_0 e^{-\gamma t}+\int_0^t e^{-\gamma(t-\xi)}\sigma e(\xi)\,d\xi
\ge i_0 e^{-\gamma t}>0
\quad\text{for every finite } t\ge 0,
\]
since $e(t)\ge 0$.

\textbf{Step 4: nonnegativity of $r$.}
Finally,
\[
r(t)=r_0+\int_0^t\big(\gamma i(\xi)+u(\xi)s(\xi)\big)\,d\xi\ge 0.
\]
\end{proof}

\begin{remark}[Role of the vaccination delay]\label{rem:delay}
The delay parameter $t_{\mathrm{delay}}$ constrains the vaccination control to $u(t)=0$ for $t < t_{\mathrm{delay}}$, reflecting the logistical lead time required for vaccine development or deployment. During this initial period, the planner must rely exclusively on the suppression instrument $h(t)$. The theoretical results of Sections~2 and~3 hold for any $t_{\mathrm{delay}} \ge 0$; the delay simply restricts the admissible set $\mathcal{U}_{\mathrm{ad}}$ without altering the structure of the existence or convergence proofs. Economically, a larger $t_{\mathrm{delay}}$ increases the shadow price of infection during the early phase, as the planner has fewer instruments available.
\end{remark}
\begin{remark}\label{Monotonicity}
Summing the first three equations of \eqref{sdeq:1} yields:
\[
(s+e+i)^{\prime}=-(\gamma i(t)+u(t) s(t)).
\]
Since $s(t)>0$ for all $t\ge 0$ and $i(t)>0$ for every finite $t\ge 0$ (by Theorem \ref{diffexist}),
the total non-recovered population is strictly decreasing.
\end{remark}

\begin{corollary}\label{cor2}
For every $t\ge 0$, the non-recovered population satisfies the integral identity:
\[
s(t)+e(t)+i(t)-(s_0+e_0+i_0)=-\int_0^t\Bigl[u(\tau)s(\tau)+\gamma\,i(\tau)\Bigr]\,d\tau.
\]
\end{corollary}

\begin{proof}
Define $X(t)=s(t)+e(t)+i(t)$. Differentiating and substituting the system dynamics yields $X'(t) = -[u(t)s(t)+\gamma i(t)]$. Integrating from 0 to $t$ yields the result.
\end{proof}

\begin{corollary}[Basic Asymptotic Properties]\label{corr3}
Under Assumption \ref{assump:bounds}, the following asymptotic properties hold:
\begin{enumerate}
    \item $s'(t) < 0$ for all $t \geq 0$, i.e., $s(t)$ is strictly decreasing.
    \item $\lim_{t \to \infty} i(t) = 0$.
    \item $\lim_{t \to \infty} e(t) = 0$.
    \item The limit $s_{\infty} := \lim_{t \to \infty} s(t)$ exists and $s_{\infty} \in [0, s_0]$.
\end{enumerate}
\end{corollary}

\begin{proof}
    \begin{enumerate}
     \item From \eqref{sdeq:1}, $s'(t) = -[(\beta - h(t))i(t) + u(t)]s(t)$. 
       We have established $s(t) > 0$ for all $t\ge 0$ and $i(t) > 0$ for every finite $t\ge 0$ (Theorem \ref{diffexist}). Also, Assumption \ref{assump:bounds} guarantees $\beta - h(t) \ge \beta - h_{\max} > 0$.
        \item From Corollary \ref{cor2}, the function $X(t) = s(t)+e(t)+i(t)$ is non-increasing and bounded below by 0. Therefore, the limit $\lim_{t \to \infty} X(t)$ exists. This implies the integral 
        \[
        \int_0^\infty [\gamma i(\tau) + u(\tau)s(\tau)] d\tau < \infty.
        \]
        Since $u, s \ge 0$, we must have $\int_0^\infty \gamma i(\tau) d\tau < \infty$, so $i \in L^1([0, \infty))$. 
        From the system dynamics, $i'(t)=\sigma e(t)-\gamma i(t)$. Since the state space $\mathbb{T}$ is compact, $e(t)$ and $i(t)$ are bounded, so $i'(t)$ is uniformly bounded. This implies $i(t)$ is uniformly continuous.
        By \textbf{Barbalat's Lemma} (a uniformly continuous function in $L^1$ with a bounded derivative must vanish at infinity), we conclude $\lim_{t \to \infty} i(t) = 0$.
        
        \item We write $e'(t) + \sigma e(t) = F(t)$, where $F(t)=(\beta-h(t))s(t)i(t)$. Note that $0\le F(t)\le \beta\, i(t)$. Since $i \in L^1$, it follows that the input $F \in L^1$. Since the kernel $e^{-\sigma t} \in L^1(\mathbb{R}_+)$ and the input $F \in L^1(\mathbb{R}_+)$, their convolution vanishes at infinity. Thus, $\lim_{t\to\infty}e(t)=0$.
        
        \item Since $s(t)$ is decreasing (Part 1) and bounded below by 0, the limit $s_\infty$ exists.
    \end{enumerate}
\end{proof}

\begin{proposition}[Final-size bound under maximal suppression]\label{prop:finalsize}
Let $\tilde\beta:=\beta-h_{\max}>0$ and $X_0:=s_0+e_0+i_0$. Under the admissible constant controls
$h(t)\equiv h_{\max}$ and $u(t)\equiv 0$, the limit $s_\infty=\lim_{t\to\infty}s(t)$ exists and satisfies
\begin{equation}\label{eq:finalsize_general}
\ln\!\Bigl(\frac{s_\infty}{s_0}\Bigr)=-\frac{\tilde\beta}{\gamma}\bigl(X_0-s_\infty\bigr).
\end{equation}
Equivalently,
\begin{equation}\label{eq:lambert_general}
s_\infty
= -\frac{\gamma}{\tilde\beta}\,
W_0\!\left(-\frac{\tilde\beta}{\gamma}\, s_0\,e^{-(\tilde\beta/\gamma)X_0}\right),
\end{equation}
where $W_0$ is the principal branch of the Lambert $W$ function.
\end{proposition}

\begin{proof}
Under $h(t)\equiv h_{\max}$ and $u(t)\equiv 0$, we have
$s'(t)=-\tilde\beta s(t)i(t)$ and $r'(t)=\gamma i(t)$. By the chain rule,
\[
\frac{ds}{dr} = \frac{ds/dt}{dr/dt} = -\frac{\tilde\beta s i}{\gamma i} = -\frac{\tilde\beta}{\gamma}s.
\]
This division is valid because $i(t) > 0$ for all finite $t$ (Theorem \ref{diffexist}).
Integrating with respect to $r$ yields $\ln(s(t)/s_0)=-(\tilde\beta/\gamma)(r(t)-r_0)$.
Letting $t\to\infty$, we use the asymptotic results $e(t)\to 0$ and $i(t)\to 0$ (Corollary~\ref{corr3}). By the conservation identity $s+e+i+r=1$, we have $X_0 = 1-r_0$ and $r_\infty = 1 - s_\infty$. Thus, $r_\infty - r_0 = (1-s_\infty) - (1-X_0) = X_0 - s_\infty$. Substituting this into the integrated equation yields \eqref{eq:finalsize_general} (see \cite{arino-2007} for standard final size derivations). Solving for $s_\infty$ yields \eqref{eq:lambert_general}.
\end{proof}

\begin{corollary}
For a general time-varying vaccination strategy $u(t)$, the final susceptible size is bounded by:
\[
s_{\infty} \leq s_0 \exp \left(-\int_0^{\infty}\left[\left(\beta-h_{\max }\right) i(\tau)+u(\tau)\right] d \tau\right).
\]
\end{corollary}
\begin{proof}
Consider the equation for $s'(t)$. Since $h(t) \le h_{\max}$ implies $-(\beta-h(t)) \le -(\beta-h_{\max})$, we have:
\[
\frac{s'(t)}{s(t)} = -(\beta-h(t))i(t) - u(t) \le -(\beta-h_{\max})i(t) - u(t).
\]
Integrating from 0 to \(\infty\) yields the inequality.
\end{proof}

\begin{lemma}[Time-Free Representation]\label{lemma:timefree}
Assume $\beta > h_{\max}$ and $i_0 > 0$. For any time $t \ge 0$, the state variables satisfy the following identity, where the explicit dependence on time is replaced by dependence on the susceptible fraction $s$:
\begin{equation}\label{eq:timefree}
(s_0+e_0+i_0) - (s(t)+e(t)+i(t)) = \int_{s(t)}^{s_0} \frac{u(\tau(\Xi)) \Xi+\gamma i(\tau(\Xi))}{\Xi[(\beta-h(\tau(\Xi))) i(\tau(\Xi))+u(\tau(\Xi))]} \, d\Xi.
\end{equation}
The expression is useful for characterizing the geometry of feasible and safe sets in Section \ref{sec:basic_results}.
\end{lemma}

\begin{proof}
The proof relies on the strict monotonicity of $s(t)$ to establish a valid change of variables.

\textbf{Step 1: Strict Monotonicity.}
Recall the state equation for $s(t)$:
\[
s'(t) = -[(\beta - h(t))i(t) + u(t)]s(t).
\]
From Theorem \ref{diffexist}, we have $s(t) > 0$ and $i(t) > 0$ for all $t \ge 0$. 
By Assumption \ref{assump:bounds}, $\beta - h(t) > 0$. Since $u(t) \ge 0$, the term in the brackets is strictly positive:
\[
(\beta - h(t))i(t) + u(t) > 0 \quad \forall t \ge 0.
\]
Consequently, $s'(t) < 0$ strictly for all $t \ge 0$. Thus, $s(t)$ is strictly decreasing and absolutely continuous on $[0, t]$, hence admits an a.e.-defined inverse $\tau$ on $[s(t), s_0]$, and the change of variables is valid by Hajłasz \cite{hajasz-1993}.

\textbf{Step 2: Change of Variables.}
We start with the integral identity derived in Corollary \ref{cor2}:
\[
\Delta X(t) := (s_0+e_0+i_0) - (s(t)+e(t)+i(t)) = \int_0^t \left[u(\tau)s(\tau)+\gamma i(\tau)\right] d\tau.
\]
We apply the change of variables $\Xi = s(\tau)$. Since $s(\tau)$ is a strictly decreasing absolute continuous function, we have $d\tau = \frac{d\Xi}{s'(\tau(\Xi))}$ almost everywhere, and the limits of integration transform from $0 \to s_0$ and $t \to s(t)$.
Substituting into the integral:
\[
\int_0^t [u(\tau)s(\tau)+\gamma i(\tau)] d\tau = \int_{s_0}^{s(t)} \frac{u(\tau(\Xi))\Xi + \gamma i(\tau(\Xi))}{s'(\tau(\Xi))} d\Xi.
\]
Substituting $s'(\tau(\Xi)) = -[(\beta-h(\tau(\Xi)))i(\tau(\Xi)) + u(\tau(\Xi))]\Xi$ and swapping the limits of integration to absorb the negative sign yields the result.
\end{proof}

\section{Dynamic Optimization}
\label{sec:dynamic_optimization}

\subsection{Contribution and Framework}
\label{subsec:contribution_framework}

Having established the well-posedness, invariance, and asymptotic properties of the controlled SEIR dynamics in the previous section, we now address the planner’s intertemporal optimization problem. A central methodological challenge in the existing epidemic-control literature is the frequent conflation of hard state constraints, such as non-negotiable ICU capacity limits, with soft penalties that merely discourage violations through cost terms. These two formulations are not mathematically equivalent as they induce different admissible sets and distinct variational objects, thereby leading to different optimality conditions and economic interpretations.
\par

The choice of an infinite planning horizon $I=[0,\infty)$ is motivated by both mathematical and practical considerations. From a practical standpoint, the terminal date of an epidemic is not known \textit{a priori}; a finite-horizon formulation introduces an artificial terminal boundary that distorts the optimal policy near the end of the planning period through transversality effects. The infinite-horizon formulation, combined with discounting at rate $\delta > 0$, eliminates this artifact and yields time-consistent policies whose structure does not depend on an arbitrarily chosen endpoint. This modeling choice is standard in the economics literature on dynamic optimization and resource management \cite{kamien-schwartz-1991, acemoglu-2009, sethi-thompson-2000}, and has been adopted in recent infinite-horizon epidemic control studies \cite{freddi-goreac-2024, della-rossa-2025, bliman-2021}. In the epidemiological context, the infinite horizon is best understood not as a claim that the epidemic literally persists forever---indeed, Corollary~\ref{corr3} guarantees $i(t)\to 0$---but rather as a device ensuring that the optimal policy accounts for all future consequences without truncation bias.
\par
This section resolves this ambiguity through a variational framework that explicitly links the two formulations. We begin by establishing the existence of optimal controls for a penalized infinite-horizon problem with discounting and bounded controls. Crucially, we introduce a theoretical bridge based on Yosida regularization, demonstrating that the quadratic over-capacity penalty term is the Moreau--Yosida envelope of the hard constraint indicator. This result proves that the optimal solutions of the penalized problem converge to the hard-constrained optimum as the penalty parameter approaches infinity \cite{RockafellarWets1998}. Furthermore, we provide an explicit topological proof of $\Gamma$-convergence to justify the use of finite-horizon approximations for the infinite-horizon problem \cite{DalMaso1993}. Finally, we characterize the optimality conditions via the Pontryagin Maximum Principle (PMP) applied to the penalized problem. We identify boundary-maintenance arcs where the switching function vanishes, showing that the planner optimally maintains the epidemic near a dynamic boundary rather than applying purely bang--bang controls. Throughout the analysis, we utilize the controlled SEIR dynamics (System~\ref{sdeq:1}) under Assumption~\ref{assump:bounds}, working with the reduced state vector $x(t):=(s(t),e(t),i(t))$ since $r(t)$ is redundant by population conservation.

\subsection{Problem Formulations}
\label{subsec:problem_formulations}

Let $I=[0,\infty)$. The reduced state is $x(t)=(s(t),e(t),i(t))\in[0,1]^3$, and the control is $v(t)=(u(t),h(t))\in[0,u_{\max}]\times[0,h_{\max}]$.
The control $u(t)$ represents the \textbf{vaccination rate} (moving mass from $s$ to $r$), while $h(t)$ represents \textbf{suppression effort} (reducing transmission from $\beta$ to $(\beta-h(t))$). Distinct from both is the \textbf{healthcare capacity constraint} $i(t)\le I_{\max}$.

\subsubsection{Admissible controls and Dynamics}
We define the admissible control set as:
\[
\mathcal{U}_{\mathrm{ad}}
:=\left\{(u,h)\in L^\infty(I;[0,u_{\max}])\times L^\infty(I;[0,h_{\max}]) \;:\;
u(t)=0 \ \text{for a.e. } t<t_{\mathrm{delay}}\right\}.
\]
Given $(u,h)\in\mathcal{U}_{\mathrm{ad}}$ and initial condition $x_0$, Theorem~\ref{diffexist} implies a unique global trajectory $x(\cdot;u,h)\in W^{1,\infty}_{\mathrm{loc}}(I;\mathbb{R}^3)$.

\subsubsection{Running cost}
Fix a discount rate $\delta>0$ and nonnegative cost coefficients $c_H,c_{NH},c_V\ge 0$. We define the base running cost density (undiscounted) as:
\[
L_0(x,v) := c_H\,i\,h+c_{NH}\,i+c_V\,u\,s.
\]
The choice of functional is grounded in specific economic reasoning. The term $c_Vus$ captures the operational cost of the vaccination program. The term $c_Hih$ captures the prevalence-scaled burden of suppression (e.g., contact tracing, enforcement resources); we interpret $c_H$ as the marginal resource cost per unit reduction in transmission intensity, scaled by current prevalence. The term $c_{NH}i$ captures the residual societal losses (morbidity, lost productivity) from infections.

\subsubsection{Hard-constrained Formulation}
Let $I_{\max}\in(0,1)$ be the capacity threshold and define the feasible state set $C:=\{x\in[0,1]^3:\ i\le I_{\max}\}$. To avoid ambiguity regarding the integration of indicator functions, we explicitly define the set of feasible controls $\mathcal{F}$ and the strict functional $J_{\mathrm{strict}}$ as follows:
\[
\mathcal{F} := \left\{ (u,h) \in \mathcal{U}_{\mathrm{ad}} : x(t; u,h) \in C \ \forall t \ge 0 \right\}.
\]
\begin{equation}\label{Pstrict}
\boxed{
J_{\mathrm{strict}}(u,h) := 
\begin{cases} 
\int_0^\infty L_0(x(t), v(t))e^{-\delta t}\,dt & \text{if } (u,h) \in \mathcal{F}, \\
+\infty & \text{otherwise.}
\end{cases}
}
\end{equation}
The hard-constrained optimization problem is $(\mathcal{P}_{\mathrm{strict}}) : \inf_{(u,h) \in \mathcal{U}_{\mathrm{ad}}} J_{\mathrm{strict}}(u,h)$.

\subsubsection{Penalized Formulation}
To facilitate numerical methods, we define the capacity-violation penalty $\psi(i):=(\max\{0,i-I_{\max}\})^2$. For a penalty weight $\kappa>0$, we define the penalized running cost density including the discount factor:
\begin{equation}\label{eq:running_cost_kappa}
\ell_\kappa(t,x,v)
:=\big(L_0(x,v) + \kappa\,\psi(i)\big)e^{-\delta t}.
\end{equation}
The penalized infinite-horizon problem is:
\begin{equation}\label{Pkappa}
\boxed{
\begin{aligned}
(\mathcal{P}_{\kappa})\qquad
& \inf_{(u,h)\in\mathcal{U}_{\mathrm{ad}}} && J_\kappa(u,h) := \int_0^\infty \ell_\kappa\bigl(t,x(t),v(t)\bigr)\,dt \\
& \text{subject to} && x(\cdot)=x(\cdot;u,h).
\end{aligned}}
\end{equation}

\subsection{Existence and Variational Convergence}
\label{subsec:existence_yosida}

\subsubsection{Topology and Convergence Framework}
We endow the admissible control set $\mathcal{U}_{\mathrm{ad}} \subset L^\infty(I; \mathbb{R}^2)$ with the weak-$\ast$ topology (see \cite{bell-1997} for discussion on optimal control in infinite horizons). We define the associated notion of convergence for sequences of controls and trajectories.

\begin{lemma}[Weak-$\ast$ to State Convergence] \label{lem:convergence}
Let $(u_n,h_n)\rightharpoonup^\ast (u,h)$ in the weak-$\ast$ topology of $L^\infty(I;\mathbb{R}^2)$ with $(u_n,h_n),(u,h)\in\mathcal{U}_{\mathrm{ad}}$. Let $x_n=x(\cdot;u_n,h_n)$ be the associated trajectories. Then there exists a subsequence (not relabeled) such that $x_n\to x(\cdot;u,h)$ uniformly on $[0,T]$ for every $T>0$.
\end{lemma}

\begin{proof}
Fix $T>0$. Since $x_n(t)\in[0,1]^3$ for all $t\ge 0$ and the vector field is bounded on $[0,1]^3\times[0,u_{\max}]\times[0,h_{\max}]$, the family $\{x_n\}$ is uniformly Lipschitz on $[0,T]$. By Arzel\`a--Ascoli, there exists a subsequence (not relabeled) and a limit $x$ such that $x_n\to x$ uniformly on $[0,T]$.

It remains to identify $x$ as the solution associated with $(u,h)$. Write the integral form for the $s$-equation on $[0,T]$:
\[
s_n(t)=s_0-\int_0^t \Big((\beta-h_n(\xi))s_n(\xi)i_n(\xi)+u_n(\xi)s_n(\xi)\Big)\,d\xi.
\]
Since $s_n\to s$ and $i_n\to i$ uniformly on $[0,T]$, we have $s_n i_n\to s i$ uniformly on $[0,T]$, hence in $L^1(0,T)$, and also $s_n\to s$ in $L^1(0,T)$. Using the decomposition
\[
\int_0^t u_n s_n\,d\xi - \int_0^t u s\,d\xi = \int_0^t u_n(s_n-s)\,d\xi + \int_0^t (u_n-u)s\,d\xi,
\]
the first term converges to $0$ because $\|u_n\|_\infty$ is uniformly bounded and $s_n\to s$ in $L^1$, while the second term converges to $0$ by weak-$\ast$ convergence of $u_n$ tested against the fixed $L^1$ function $s\,\mathbf{1}_{[0,t]}$ (where $\mathbf{1}_{[0,t]}(\xi) = 1$ if $\xi \in [0,t]$). Similarly,
\[
\int_0^t h_n(s_n i_n)\,d\xi - \int_0^t h(s i)\,d\xi = \int_0^t h_n(s_n i_n - s i)\,d\xi + \int_0^t (h_n-h)(s i)\,d\xi,
\]
where the first term vanishes by boundedness of $h_n$ and $L^1$ convergence of $s_n i_n$, and the second term vanishes by weak-$\ast$ convergence of $h_n$ tested against the fixed $L^1$ function $(s i)\mathbf{1}_{[0,t]}$. Therefore the integral equation for $s_n$ passes to the limit, and the limit $s$ satisfies the integral equation with controls $(u,h)$. The same argument applies to the $e$ and $i$ equations, hence $x$ coincides with $x(\cdot;u,h)$ on $[0,T]$.
\end{proof}

\begin{theorem}[Existence of an optimal solution for $\mathcal{P}_\kappa$]\label{thm:exist_kappa}
Under Assumption~\ref{assump:bounds}, for every $\kappa>0$ there exists an optimal control $(u^\kappa,h^\kappa)\in\mathcal{U}_{\mathrm{ad}}$ solving \eqref{Pkappa}.
\end{theorem}

\begin{proof}
Let $\{(u_n,h_n)\}_{n \in \mathbb{N}} \subset \mathcal{U}_{\mathrm{ad}}$ be a minimizing sequence for $J_\kappa$.

\begin{enumerate}
    \item \textbf{Compactness:} The set $\mathcal{U}_{\mathrm{ad}}$ is a subset of the ball $\{v\in L^\infty(I;\mathbb R^2):\|v\|_\infty\le M\}$ where $M=\max\{u_{\max}, h_{\max}\}$. This ball is weak-* compact by Banach--Alaoglu. Since $L^1(I)$ is separable, the weak-* topology on this ball is metrizable, implying weak-* sequential compactness. The box constraints and the linear delay constraint are weak-* closed (for the delay: for any test function $\varphi \in L^1(I)$ supported on $[0, t_{\text{delay}})$, we have $\int u_n \varphi\, dt = 0$ for all $n$, which implies $\int u \varphi\, dt = 0$ by weak-$\ast$ convergence, so $u=0$ a.e.\ on the delay interval). Thus, $\mathcal{U}_{\mathrm{ad}}$ is weak-* sequentially compact. There exists a subsequence (not relabeled) such that $(u_n,h_n)\rightharpoonup^\ast(u^\kappa,h^\kappa) \in \mathcal{U}_{\mathrm{ad}}$.
    
    \item \textbf{Convergence of Trajectories:} By Lemma~\ref{lem:convergence}, we extract a further subsequence such that the corresponding state trajectories $x_n$ converge uniformly on compact sets to $x^\kappa = x(\cdot; u^\kappa, h^\kappa)$.
    
    \item \textbf{Term-by-Term Convergence:} We show $J_\kappa(u^\kappa,h^\kappa)\le \liminf_n J_\kappa(u_n,h_n)$. The integrand involves bilinear terms ($u s$, $i h$) and the penalty. For any finite $\bar{T} > 0$, consider the term $\int_0^{\bar{T}} u_n(t) s_n(t) e^{-\delta t}\, dt$. We decompose it as:
    \[
    \int_0^{\bar{T}} u_n(s_n - s^\kappa)e^{-\delta t}\, dt + \int_0^{\bar{T}} (u_n - u^\kappa)s^\kappa e^{-\delta t}\, dt + \int_0^{\bar{T}} u^\kappa s^\kappa e^{-\delta t}\, dt.
    \]
    The first integral vanishes because $u_n$ is bounded and $s_n \to s^\kappa$ uniformly. The second vanishes because $(u_n - u^\kappa) \rightharpoonup^\ast 0$ and $s^\kappa e^{-\delta t} \in L^1([0,\bar{T}])$. The same decomposition applies to $i_n h_n$ and $i_n$. Additionally, $\psi(i_n) \to \psi(i^\kappa)$ uniformly on $[0,\bar{T}]$, ensuring convergence of the penalty integral. Thus, the integrals converge on $[0,\bar{T}]$. Finally, since $i(t) \in [0,1]$, we have $0 \le \psi(i(t)) \le (1-I_{\max})^2$, hence the integrand is bounded by $C_\kappa e^{-\delta t} \in L^1([0,\infty))$, and the tail contribution vanishes as $\bar{T} \to \infty$. (See also \cite{dacorogna-2007} for general existence results via direct methods).
\end{enumerate}
Therefore, $(u^\kappa,h^\kappa)$ is an optimal control.
\end{proof}

\subsubsection{Moreau--Yosida regularization: linking soft and hard constraints}
We now rigorously link the two problems. The hard constraint set $C=\{x\in[0,1]^3:\ i\le I_{\max}\}$, as defined in Section~\ref{subsec:problem_formulations}, can be represented by the extended-value indicator function $\iota_C(x)$. Pointwise in the scalar variable $i$, the Moreau--Yosida envelope of the scalar indicator $\iota_{(-\infty,I_{\max}]}(i)$ with parameter $\varepsilon>0$ is defined as:
\[
\inf_{y\le I_{\max}}\left\{\iota_{(-\infty,I_{\max}]}(y)+\frac{1}{2\varepsilon}(i-y)^2\right\} = \frac{1}{2\varepsilon}\big(\max\{0,i-I_{\max}\}\big)^2.
\]
Our penalty term $\kappa\psi(i)$ corresponds exactly to this envelope with parameter $\varepsilon=(2\kappa)^{-1}$. This implies that as $\kappa\to\infty$, the penalty term epi-converges to the indicator function of the hard constraint.

\begin{theorem}[Moreau--Yosida / $\Gamma$-limit as $\kappa\to\infty$]\label{thm:yosida_limit}
Assume $\mathcal{P}_{\mathrm{strict}}$ admits at least one feasible pair. Let $\{\kappa_n\}$ be any sequence with \(\kappa_n \to \infty\). Then \(J_{\kappa_n}\) \(\Gamma\)-converges to \(J_{\mathrm{strict}}\) in the weak-\(\ast\) topology on \(\mathcal{U}_{\mathrm{ad}}\). Consequently, any sequence of minimizers of \(J_{\kappa_n}\) admits a subsequence converging to a minimizer of \(J_{\mathrm{strict}}\).
\end{theorem}

\begin{proof}
Since $\mathcal{U}_{\mathrm{ad}}$ is weak-$\ast$ compact, the functionals are equi-coercive. We verify the two $\Gamma$-convergence conditions:

\textbf{1. Liminf inequality:} Let $(u_n,h_n) \rightharpoonup^\ast (u,h)$. Choose a subsequence $\{n_j\}$ such that $J_{\kappa_{n_j}}(u_{n_j}, h_{n_j}) \to \liminf_{n \to \infty} J_{\kappa_n}(u_n, h_n)$. Let $x_{n_j} \to x$ uniformly on compacts (by Lemma~\ref{lem:convergence}). If the limit trajectory $x(\cdot)$ violates the constraint, there exists $t_0$ such that $i(t_0) > I_{\max}$. By continuity, there exist $\eta, \varepsilon > 0$ such that $i(t) \ge I_{\max} + \eta$ on $[t_0-\varepsilon, t_0+\varepsilon]$. Uniform convergence implies that for $n_j$ large enough, $i_{n_j}(t) \ge I_{\max} + \eta/2$ on this interval. Thus:
\[
\int_0^\infty \kappa_{n_j} \psi(i_{n_j}(t))e^{-\delta t} \, dt \ge \kappa_{n_j} \left(\frac{\eta}{2}\right)^2 \int_{t_0 - \varepsilon}^{t_0 + \varepsilon} e^{-\delta t} \, dt \xrightarrow{j \to \infty} +\infty,
\]
consistent with $J_{\mathrm{strict}} = \infty$. If $i(t) \le I_{\max}$ everywhere, the penalty vanishes, and the inequality follows from the same term-by-term decomposition used in Theorem \ref{thm:exist_kappa}, applied on $[0,M]$ and then letting $M\to\infty$.

\textbf{2. Recovery sequence:}Let $(u,h) \in \mathcal{F}$ (the set of feasible controls defined in Section~\ref{subsec:problem_formulations}). Consider the constant sequence $(u_n,h_n)=(u,h)$. Since the trajectory is feasible, $\psi(i(t;u,h))\equiv 0$. Thus, $J_{\kappa_n}(u, h) = J_{\mathrm{strict}}(u,h)$ for all $n$.
By the Fundamental Theorem of $\Gamma$-convergence \cite{DalMaso1993}, minimizers converge.
\end{proof}

\subsection{Variational convergence: finite to infinite horizon}\label{subsec:gamma_T}

To justify the large-time horizon numerical approximations, we define the finite-horizon problem for each $T>0$:
\[
(\mathcal{P}_{\kappa,T}): \quad \inf_{(u,h) \in \mathcal{U}_{\mathrm{ad}}} J_{\kappa,T}(u,h):=\int_0^T \ell_\kappa\bigl(t,x(t;u,h),v(t)\bigr)\,dt.
\]

\begin{theorem}[$\Gamma$-convergence as $T\to\infty$]\label{thm:gamma_T}
Fix $\kappa>0$ and $\delta>0$. Let $\{T_n\}$ be any sequence with $T_n \to \infty$. Then $J_{\kappa,T_n}$ $\Gamma$-converges to $J_\kappa$ in the weak-$\ast$ topology on \(\mathcal{U}_{\mathrm{ad}}\). Moreover, any sequence of minimizers of \((\mathcal{P}_{\kappa,T_n})\) admits a subsequence converging to a minimizer of \(\mathcal{P}_\kappa\).
\end{theorem}

\begin{proof}

\textbf{1. Liminf inequality:} Let $(u_n,h_n) \rightharpoonup^\ast (u,h)$. Choose a subsequence $\{n_j\}$ attaining the liminf. Along a further subsequence, $x_{n_j} \to x$ uniformly on compacts. For any fixed $M>0$, eventually $T_{n_j} \ge M$, so $J_{\kappa,T_{n_j}}(u_{n_j}, h_{n_j}) \ge \int_0^M \ell_\kappa(t, x_{n_j}, v_{n_j}) dt$. Taking the limit as $j \to \infty$ (implying uniform convergence on $[0,M]$), and then letting $M \to \infty$, yields the inequality.

\textbf{2. Recovery sequence:} Take the constant sequence $u_n=u, h_n=h$. The difference $|J_\kappa - J_{\kappa,T_n}|$ is the tail integral $\int_{T_n}^\infty \ell_\kappa\bigl(t,x(t;u,h),v(t)\bigr)\,dt$. Since the state space is compact and controls are bounded, $\ell_\kappa(t,x,v) \le C e^{-\delta t} \in L^1$. By the Dominated Convergence Theorem, the tail integral vanishes as $n\to\infty$.
\end{proof}

\subsection{Optimality conditions and economic interpretation}\label{subsec:pmp_econ}

We apply the Pontryagin Maximum Principle (PMP) to the finite-horizon penalized problem $(\mathcal{P}_{\kappa,T})$ on $[0,T]$, with free terminal state (implying the transversality condition $\lambda(T)=0$). Define the Hamiltonian:
\[
H_\kappa(t,x,v,\lambda)
= \ell_\kappa(t,x,v) + \lambda \cdot f(x,v)
= \ell_\kappa(t,x,v) + \lambda_s \dot{s} + \lambda_e \dot{e} + \lambda_i \dot{i}.
\]
Since \(\psi(i)=(\max\{0,\,i-I_{\max}\})^2\) is continuously differentiable, the adjoint equations hold in the classical sense (see \cite{mangasarian-1966} for sufficient conditions, \cite{avram-2021} for PMP in similar contexts, and \cite{aseev-kryazhimskiy-2007, carlson-haurie-leizarowitz-1991, aseev-kryazhimskiy-2004} for comprehensive treatments of the Pontryagin Maximum Principle on infinite-horizon problems).

\subsubsection{Adjoint system and Shadow Prices}
The adjoint equations $\dot{\lambda} = - \nabla_x H_\kappa$ are:
\begin{align}
\dot \lambda_s &= -c_V u e^{-\delta t} + \lambda_s[(\beta-h)i + u] - \lambda_e(\beta-h)i, \label{adj_s}\\
\dot \lambda_e &= \sigma(\lambda_e-\lambda_i), \label{adj_e}\\
\dot \lambda_i &= -\left[ c_H h + c_{NH} + 2\kappa\max\{0,\,i-I_{\max}\} \right] e^{-\delta t} + (\lambda_s-\lambda_e)(\beta-h)s + \gamma \lambda_i. \label{adj_i}
\end{align}
Economically, $\lambda_s(t)$ represents the \textbf{marginal shadow cost of susceptibility} (the reduction in future objective value achievable by removing an individual from the susceptible compartment), while \(\lambda_i(t)\) represents the \textbf{marginal social cost of infection}. The term $2\kappa\max\{0,\,i-I_{\max}\}$ acts as a soft barrier, signaling high costs when capacity is breached.

\subsubsection{Switching functions and Optimality}
The optimal controls $(u^*, h^*)$ minimize the Hamiltonian $H_\kappa$ pointwise. Since controls enter linearly, we define the switching functions $\Phi_u = \partial H_\kappa / \partial u$ and \(\Phi_h = \partial H_\kappa / \partial h\). Note that the term involving \(h\) in the dynamics contributes \(\lambda_s(si) + \lambda_e(-si) = si(\lambda_s - \lambda_e)\). Since \(H_\kappa\) is affine in \(u\) and \(h\), the minimizer is attained at the bounds depending on the sign of \(\Phi\):
\begin{align}
\Phi_u(t) &:= s(t)\big(c_V e^{-\delta t}-\lambda_s(t)\big), \\
\Phi_h(t) &:= i(t)\Big(c_H e^{-\delta t} + s(t)\big(\lambda_s(t)-\lambda_e(t)\big)\Big).
\end{align}
The optimal policy follows the bang-bang law:
\[
u^*(t) \in \begin{cases} \{u_{\max}\} & \Phi_u(t) < 0, \\ [0, u_{\max}] & \Phi_u(t) = 0, \\ \{0\} & \Phi_u(t) > 0. \end{cases}
\qquad
h^*(t) \in \begin{cases} \{h_{\max}\} & \Phi_h(t) < 0, \\ [0, h_{\max}] & \Phi_h(t) = 0, \\ \{0\} & \Phi_h(t) > 0. \end{cases}
\]

\subsubsection{Boundary-Maintenance Arcs}
In optimal control theory, a \emph{singular arc} refers to an interval on which a switching function vanishes identically, so that the optimal control is not determined by the bang-bang law alone. A distinct but related phenomenon is the \emph{boundary-maintenance arc}, which arises when the optimal trajectory reaches a state (phase) constraint and is held at its boundary. We analyze this latter regime, in which the infection level is maintained at the capacity limit.
\[
i(t) \equiv I_{\max} \quad \text{on some interval}.
\]
Then $\dot i(t)=0$ on that interval, hence the dynamics imply the boundary identity $\sigma e(t) = \gamma I_{\max}$. In applications it is natural to further consider a quasi-steady boundary regime in which the exposed class does not drift systematically along the boundary, i.e.\ $\dot e(t)=0$ in addition to $i(t)\equiv I_{\max}$. Under this additional restriction, substituting $\sigma e = \gamma I_{\max}$ and $\dot e = 0$ into the $e$-equation yields the explicit feedback law for suppression:
\[
h_{\mathrm{bm}}(t)=\beta-\frac{\gamma}{s(t)} \quad \iff \quad \mathcal{R}_{\mathrm{eff}}(t) := \frac{\bigl(\beta-h_{\mathrm{bm}}(t)\bigr)}{\gamma}s(t)=1.
\]
Here $\mathcal{R}_{\mathrm{eff}}(t)$ is the \textbf{effective reproduction number} at time $t$: the expected number of secondary infections produced by a single infected individual given the current susceptible fraction $s(t)$ and suppression level $h_{\mathrm{bm}}(t)$. The condition $\mathcal{R}_{\mathrm{eff}}(t)=1$ characterizes the boundary-maintenance arc as the regime in which the epidemic is held in exact equilibrium, neither growing nor declining.

This arc is admissible only when $0\le h_{\mathrm{bm}}(t)\le h_{\max}$, equivalently when
\[
\frac{\gamma}{\beta}\le s(t)\le \frac{\gamma}{\beta-h_{\max}}.
\]
Outside this range, the boundary cannot be maintained by admissible suppression alone, and the constraint forces a corner regime.

\subsection{Analytical sensitivity: marginal value of capacity}\label{subsec:sensitivity_capacity}
We define the value function $V_\kappa(u_{\max}, h_{\max})$ as the infimum of the penalized functional $J_\kappa(u, h)$ subject to the dynamic constraints and the control bounds $0 \le u(t) \le u_{\max}$ and $0 \le h(t) \le h_{\max}$. To quantify the economic value of relaxation in the infrastructure constraints (i.e., increasing vaccination supply $u_{\max}$ or hospital surge capacity $h_{\max}$), we employ a formal sensitivity analysis.

Let $\nu_u(t)$ and $\nu_h(t)$ denote the non-negative Lagrange multipliers (dual variables) associated with the upper bound constraints $u(t) \le u_{\max}$ and $h(t) \le h_{\max}$, respectively. Under standard regularity assumptions and constraint qualifications that ensure the differentiability of the value function with respect to parameters, the Envelope Theorem yields:
\begin{equation}\label{eq:envelope_thm}
\frac{\partial V_\kappa}{\partial u_{\max}} = -\int_0^\infty \nu_u(t)\,dt, \qquad
\frac{\partial V_\kappa}{\partial h_{\max}} = -\int_0^\infty \nu_h(t)\,dt.
\end{equation}
These relations offer a direct economic interpretation:
\begin{itemize}
    \item The multipliers $\nu_u(t)$ and $\nu_h(t)$ measure the instantaneous shadow price of the capacity constraints. Specifically, $\nu_u(t) > 0$ only when the optimal vaccination rate is binding at the limit ($u^*(t) = u_{\max}$).
    \item The integral $\int_0^\infty \nu_u(t) dt$ aggregates these instantaneous shadow prices over the entire planning horizon.
    \item Consequently, the marginal benefit of investing in additional capacity is proportional to the cumulative duration and intensity with which the optimal intervention is constrained by the current limits. If the constraint is rarely binding (e.g., in a "safe set" scenario), the marginal value vanishes; if the constraint binds frequently (e.g., during peak infection periods), the marginal value is substantial.
\end{itemize}
This result requires standard regularity conditions (e.g., uniqueness of multipliers, differentiability of the value function) which we assume hold for the purpose of economic interpretation.

\section{Numerical Simulations}\label{sec:simulation}

This section validates the theoretical findings established in the preceding sections through a series of numerical experiments. The analysis focuses on (i) the comparative performance of heuristic intervention strategies, and (ii) the economic behavior of the system under the optimal control policy derived via the Pontryagin Maximum Principle. 

\subsection{Comparative Analysis of Intervention Strategies}
\label{subsec:comparative_analysis}

We first evaluate the epidemic's progression under four intervention regimes: no intervention,
vaccination only, suppression only, and a combined strategy. Simulations are initialized with
$s_0 = 0.90$, $e_0 = 0.05$, $i_0 = 0.05$, and $r_0 = 0.00$ representing an early-stage outbreak in a
vulnerable population. The baseline parameters are set at $\beta = 0.5$, $\sigma = 0.2$, and
$\gamma = 0.1$, with constant control rates for vaccination ($u = 0.05$) and suppression ($h = 0.2$).

\begin{figure}[H]
    \centering
    \includegraphics[width=0.95\textwidth]{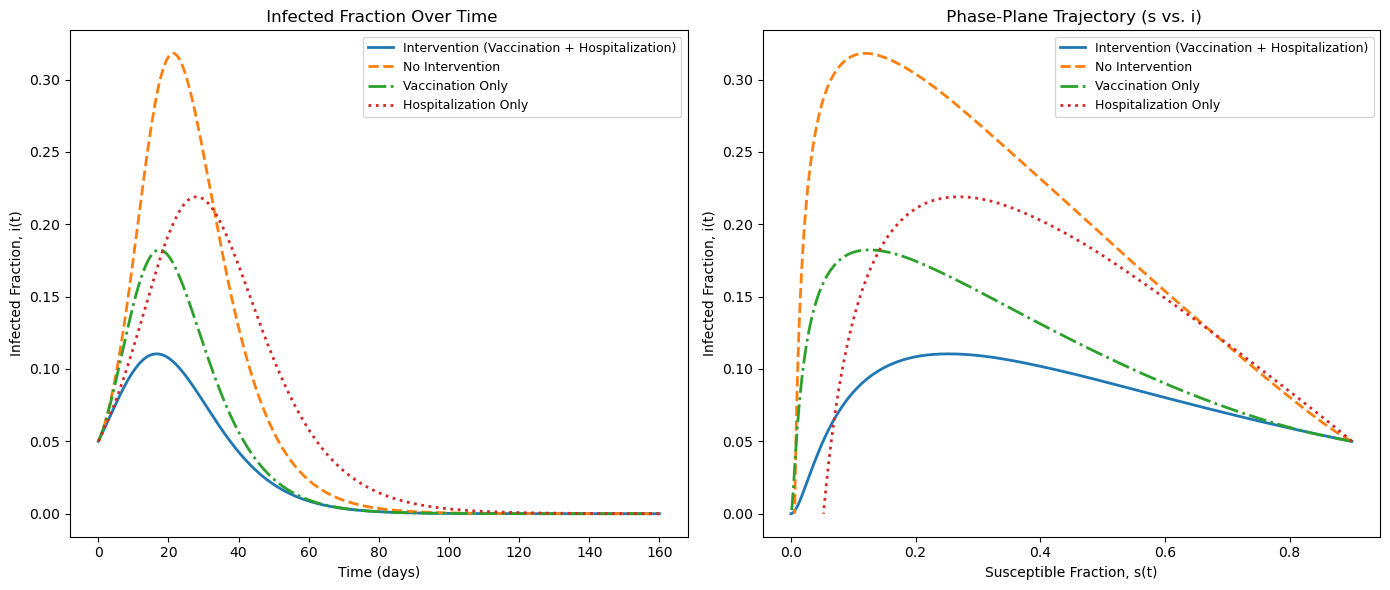}
    \caption{Epidemic dynamics under fixed intervention strategies (left) and corresponding phase-plane trajectories (right).} 
    \label{fig:interventions}
\end{figure}

As illustrated in Figure \ref{fig:interventions}, the absence of intervention results in a rapid infection surge, with $i(t)$ peaking at approximately 0.32. While individual controls provide moderate relief, the combined strategy achieves the lowest infection peak and and preserves the largest fraction of the susceptible population relative to the other scenarios.

To assess the impact of timing, we compare these results with two time-varying schedules in Table \ref{tab:peak-final}. Schedule \#1 (Strong Early) emphasizes an aggressive initial response, while Schedule \#2 (Ramp Up) delays intensity.  The results confirm that the timing of interventions is as critical as their magnitude: Schedule \#1 suppresses the peak to 0.0731, effectively keeping the outbreak within healthcare capacity limits. %

\begin{table}[H]
    \centering
    \caption{Peak Infection Fraction and Final Epidemic Size ($1 - s(\infty)$)} %
    \label{tab:peak-final}
    \begin{tabular}{lcc}
        \toprule
        \textbf{Scenario} & \textbf{Peak $i(t)$} & \textbf{Final Size} \\
        \midrule
        Constant Controls       & 0.1104 & 0.9999 \\ %
        Schedule \#1 (Strong Early)  & 0.0731 & 0.9992 \\ %
        Schedule \#2 (Ramp Up)  & 0.2045 & 0.9990 \\ %
        \bottomrule
    \end{tabular}
\end{table}

\subsection{Optimal Control Dynamics and Shadow Prices}
\label{subsec:optimal_control}

We now examine the system under the optimal control derived in Section 3. To maintain the infection fraction near the designated capacity $I_{\max} = 0.10$, the optimal vaccination rate $u(t)$ saturates at its maximum bound of 0.050 for the first 30 days.  Simultaneously, the optimal suppression $h(t)$ starts at a maximum of 0.200 before decreasing as the curve flattens. 

\begin{figure}[H]
    \centering
    \includegraphics[width=0.85\textwidth]{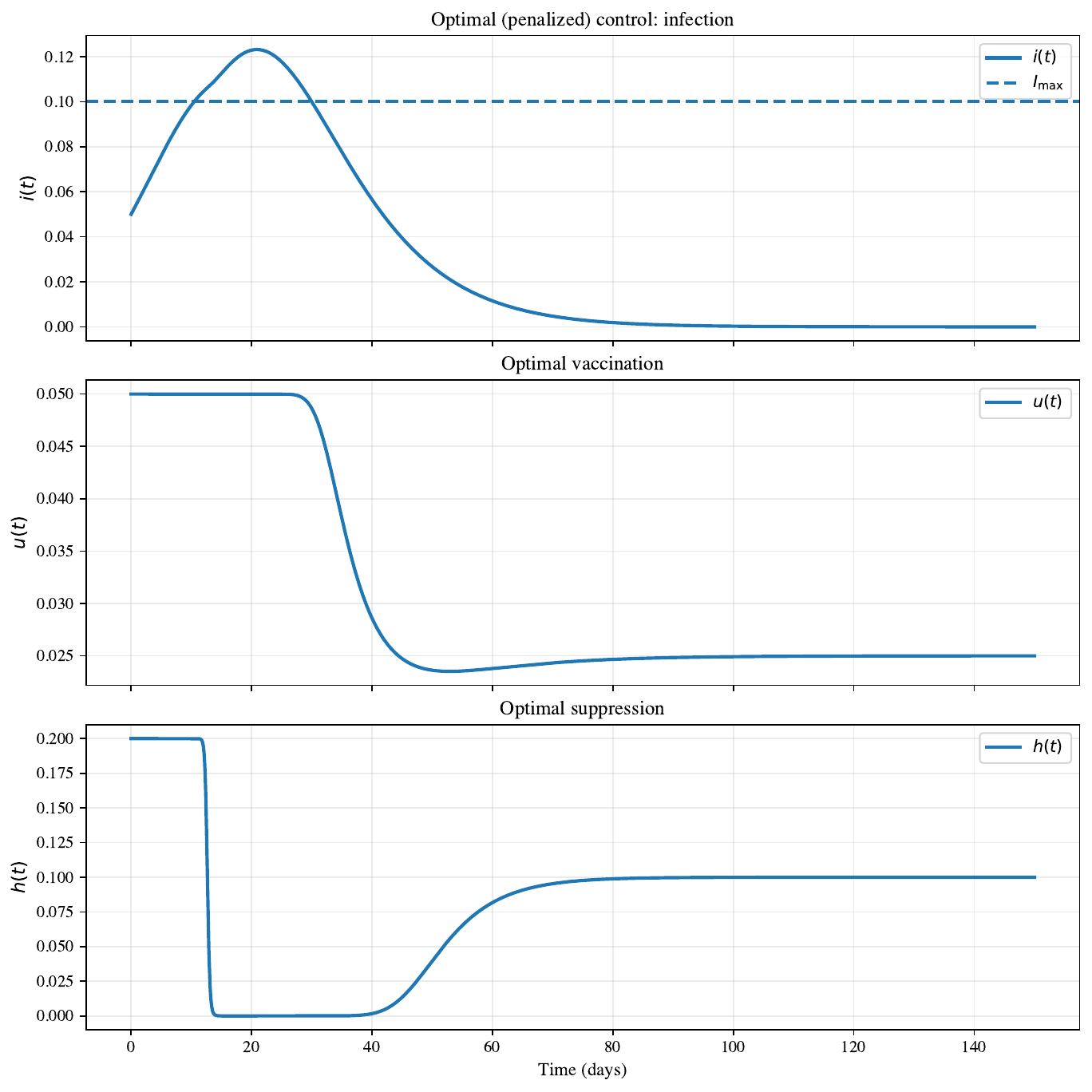}
    \caption{Optimal control trajectories: Infected fraction relative to $I_{max}$ (top), optimal vaccination $u(t)$ (middle), and optimal suppression $h(t)$ (bottom).} 
    \label{fig:optimal_controls}
\end{figure}

\begin{figure}[H]
    \centering
    \includegraphics[width=0.8\textwidth]{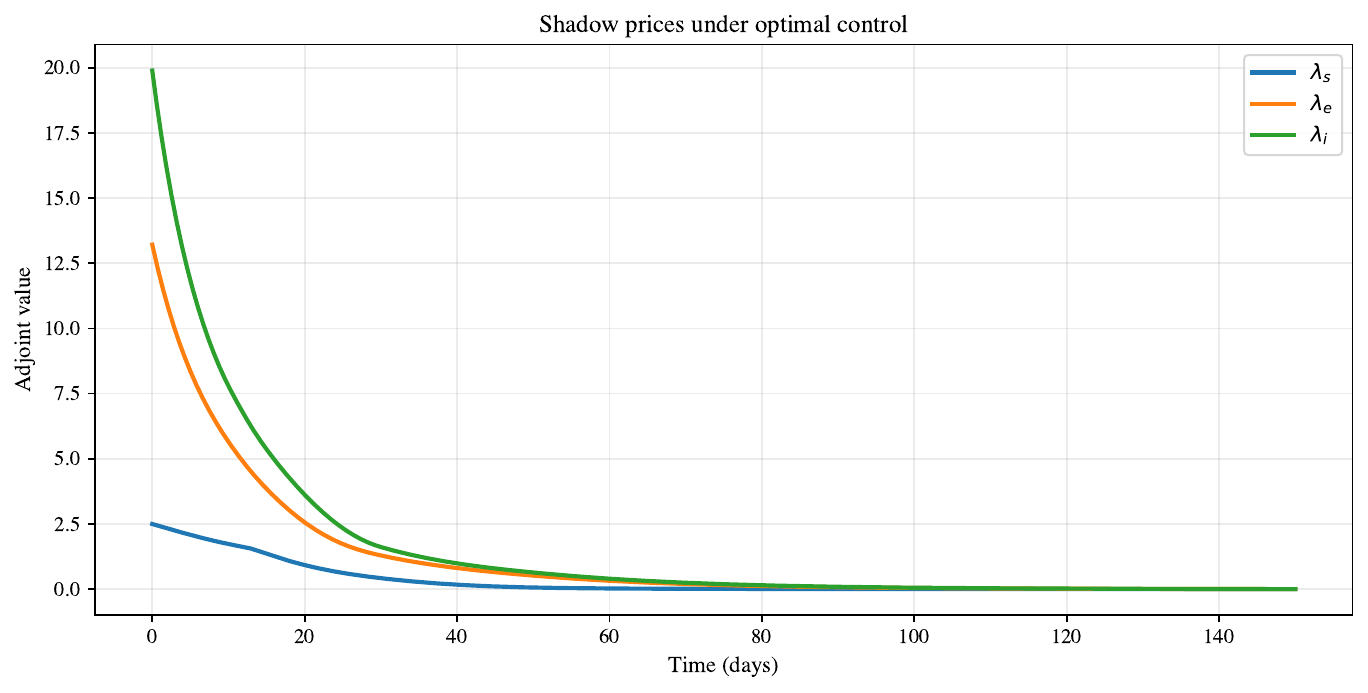}
    \caption{Time path of shadow prices $\lambda_s(t)$, $\lambda_e(t)$, and $\lambda_i(t)$ under optimal control.}
    \label{fig:shadow_prices}
\end{figure}

The economic rationale for this policy is revealed by the adjoint variables, or shadow prices, in Figure \ref{fig:shadow_prices}. The shadow price of infection, $\lambda_i(t)$, represents the marginal cost of an additional infected individual; it begins at a peak value of approximately 20.0 and decays as the epidemic is brought under control.  The trajectory of $\lambda_s(t)$ represents the marginal value of preserving the susceptible mass to avoid downstream costs. It should be noted that the optimal control simulation presented in Figure~\ref{fig:optimal_controls} uses $t_{\mathrm{delay}} = 0$, so that vaccination is available from the outset. The impact of positive delays on total cost is captured in the sensitivity analysis of Section~\ref{subsec:sensitivity} (Figure~\ref{fig:boxplots}), which demonstrates that increasing $t_{\mathrm{delay}}$ substantially raises the median total cost. A detailed comparison of optimal trajectories with and without the delay is an important direction for future numerical investigation.

\subsection{Sensitivity Analysis and Policy Implications}
\label{subsec:sensitivity}

\begin{figure}[H]
    \centering
    \includegraphics[width=0.7\textwidth]{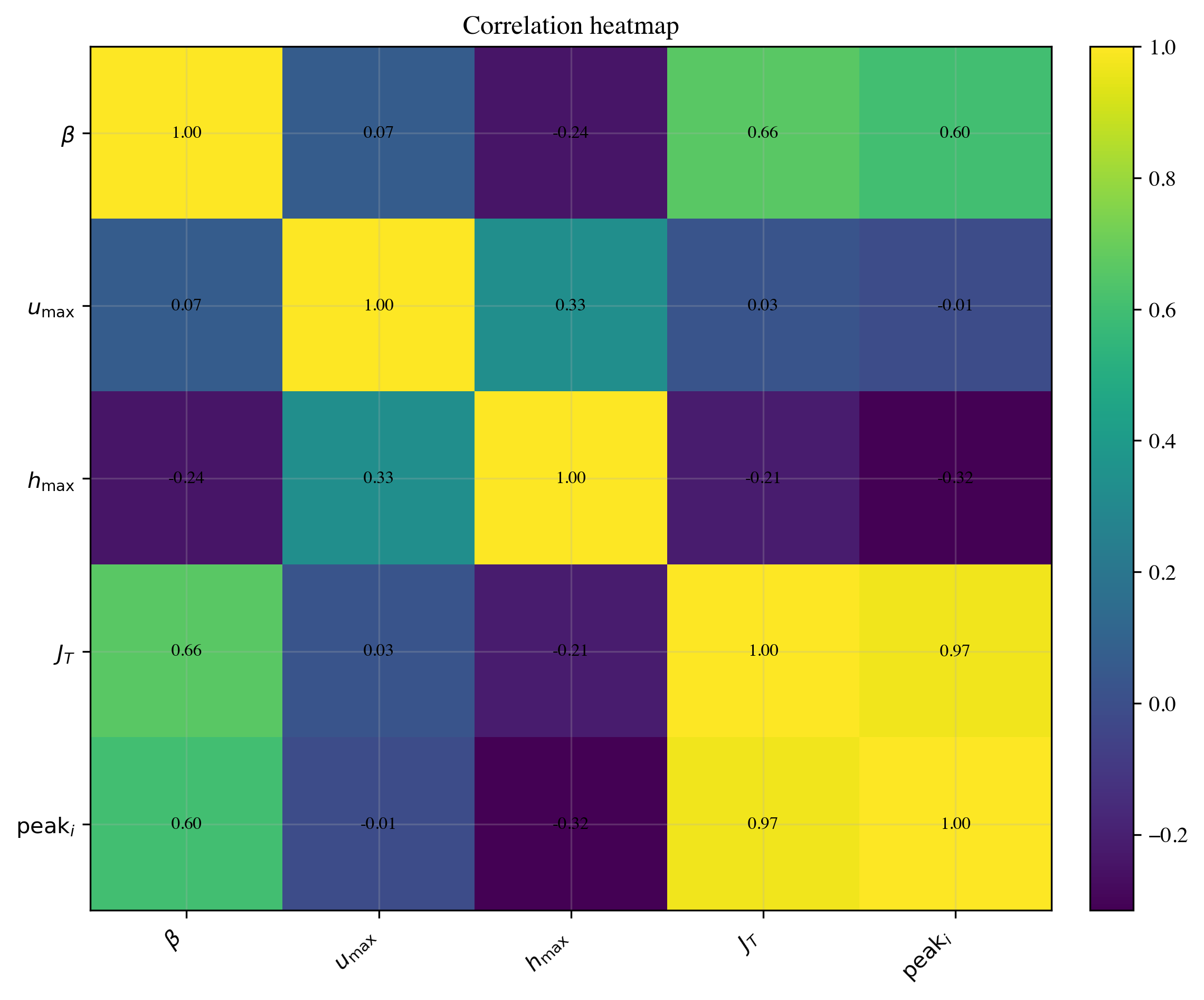}
    \caption{Correlation heatmap of parameters ($\beta, u_{max}, h_{max}$) and outcomes ($J_T, peak_i$).} 
    \label{fig:corr_heatmap}
\end{figure}

A sensitivity analysis of the total cost $J_T$ and the peak infection fraction reveals critical dependencies. As shown in the correlation heatmap (Figure \ref{fig:corr_heatmap}), the transmission rate $\beta$ has a significant positive correlation (0.66) with the total cost. Most notably, $J_T$ and the infection peak exhibit a near-perfect correlation of 0.97.

\begin{figure}[H]
    \centering
    \includegraphics[width=0.9\textwidth]{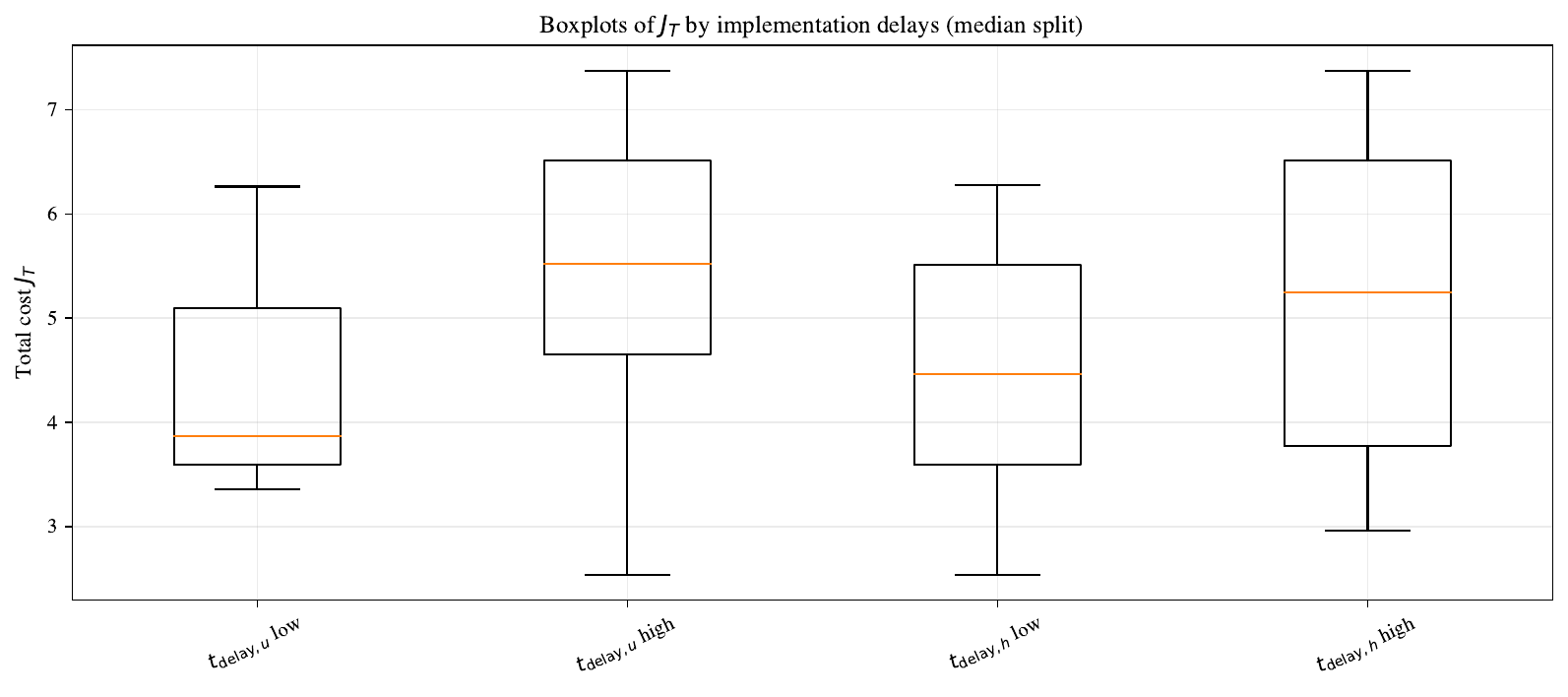}
    \caption{Boxplots of $J_T$ by implementation delays, illustrating the high cost of delayed policy execution.} 
    \label{fig:boxplots}
\end{figure}

Furthermore, the analysis of implementation delays (Figure \ref{fig:boxplots}) indicates that the cost of late intervention is non-linear. High delays in either vaccination ($t_{delay,u}$) or hospitalization ($t_{delay,h}$) significantly shift the median total cost upward. This validates the policy implication that maximizing the speed and scale of interventions provides far greater economic benefit than minimizing per-unit operational costs. \footnote{For computational details in Python, please refer to \href{https://github.com/BehroozMoosavi/Codes/blob/main/Epidemic\%20With\%20Intervention/Epidemic.ipynb}{my github}.}

\section{Conclusion and Future Directions}
\label{sec:conclusion}

This paper has established a rigorous optimal control framework for epidemic management that bridges the gap between theoretical well-posedness and operational reality. By extending the classical SEIR compartmental model to incorporate implementation delays, vaccination logistics, and hard healthcare capacity constraints, we have provided a prescriptive tool for balancing public health stability with socio-economic costs.

\par

From a theoretical perspective, our contribution is threefold. First, we proved the global existence, positivity, and asymptotic stability of the controlled delay differential system, utilizing Barbalat's Lemma to ensure the long-term convergence of the infected population. Second, we resolved the analytical difficulty of state constraints (ICU capacity) by employing Moreau--Yosida regularization. This approach allowed us to approximate the hard constraint via a sequence of penalized unconstrained problems, providing a robust justification for the numerical methods used. Third, we utilized $\Gamma$-convergence to rigorously link finite-horizon optimal control problems to the infinite-horizon setting, justifying the use of numerical approximations to infer long-term policy structures.
\par
Analytically, we characterized the optimal intervention strategy using the Pontryagin Maximum Principle, identifying boundary-maintenance arcs where the optimal policy maintains the epidemic trajectory precisely at the healthcare capacity limit. This "boundary-maintenance" regime represents a critical operational insight often missed by purely descriptive models. Furthermore, our derivation of a time-free representation of the dynamics and the use of the Lambert $W$ function provided sharp bounds on the final epidemic size under maximal suppression.

\par

Our numerical simulations illustrate these theoretical findings, quantifying the high shadow prices of infection during the peak phase and demonstrating the non-linear cost escalation associated with intervention delays. Sensitivity analysis further highlighted that expanding vaccination capacity $u_{\max}$ yields the highest marginal reduction in total social cost, emphasizing the priority of logistical speed over magnitude.

\par

Despite these contributions, this study is subject to certain limitations that suggest avenues for future research. First, the assumption of a well-mixed population simplifies the complex contact structures inherent in real-world transmission; extending this framework to metapopulation models or network-based dynamics would enhance its spatial applicability. Second, while we incorporated deterministic delays, the model assumes fixed biological parameters. Integrating stochastic differential equations to model parameter uncertainty—such as viral mutations or shocks to the transmission rate $\beta$—would better capture the unpredictability of epidemic evolution. Finally, while the focus of this work was theoretical rigor, calibrating the model against empirical data from specific outbreaks (e.g., COVID-19 or seasonal influenza) remains a necessary step to translate these optimal control protocols into actionable public health policy.

\bibliographystyle{abbrv}
\bibliography{References}

\end{document}